\documentclass[10pt]{article}
\usepackage{graphicx}
\usepackage{amsmath}
\usepackage[dvips]{epsfig}
\usepackage{amssymb}
\begin{document}
\begin{center} {\Large \bf  A family of generalized $q$-Genocchi numbers and polynomials}
\\ \vspace*{12 true pt}  T.  Kim$^1$, Byungje Lee$^2$,   C. S.  Ryoo$^3$
\vspace*{12 true pt} \\
$^1$ Division of General Education-Mathematics,
 Kwangwoon University, Seoul 139-701,  Korea \\
 $^2$ Depart. of Wireless Communications Engineering,
 Kwangwoon University, Seoul 139-701,  Korea \\
              $^3$ Department of Mathematics, Hannam University,  Daejeon 306-791, Korea \\
          \end{center}
\vspace*{12 true pt} \noindent {\bf Abstract :} In this paper we
first consider the $q$-extension of the generating function for
the higher-order generalized Genocchi numbers and polynomials
attached to $\chi$. The purpose of this paper is to present a
systemic study of some families of higher-order generalized
$q$-Genocchi numbers and polynomials attached to $\chi$ by using
the generating function of those numbers and polynomials.

\vspace*{12 true pt} \noindent {\bf 2000 Mathematics Subject
Classification :} 11B68, 11S40, 11S80

 \vspace*{12 true pt} \noindent {\bf
Key words :}   Genocchi numbers, Genocchi  polynomials,
higher-order generalized $q$-Genocchi numbers and polynomials

\begin{center} {\bf 1. Introduction } \end{center}
As well known definition, the Genocchi polynomials are defined by

$$ \left( \dfrac{2t}{ e^t +1} \right) e^{xt}= e^{G(x)t}= \sum_{n=0}^{\infty} G_{n}(x) \dfrac{t^n}{n!}, \quad | t| < \pi,  \eqno(1)$$
where we use the technical method's notation by replacing $G^n(x)$
by $G_n(x)$, symbolically, (see [3,11]). In the special case
$x=0$, $G_n=G_n(0)$ are called the $n$-th Genocchi numbers. From
the definition of Genocchi numbers, we note that $G_1=1,
G_3=G_5=G_7= \cdots =0$, and even coefficients are given by
$G_{2n}=2(1-2^{2n}) B_{2n}=2n E_{2n-1}(0)$, (see [8]), where $B_n$
is a Bernoulli number and $E_n(x)$ is an Euler polynomial. The
first few Genocchi numbers for $2, 4, 6, \cdots $ are $-1, 1, -3,
17, -155, 2073, \cdots.$ The first few prime Genocchi numbers are
given by $G_6=-3$ and $G_8=17.$ It is known that there are no
others prime Genocchi numbers with $n < 10^5$. For a real or
complex parameter $\alpha$, the higher-order Genocchi polynomials
are defined by
$$ \left( \dfrac{2t}{ e^t +1} \right)^{\alpha} e^{xt}= \sum_{n=0}^{\infty} G_{n}^{(\alpha)}(x) \dfrac{t^n}{n!}, \text{ (see [3, 4])}. \eqno(2)$$
In the special case $x=0$, $G_n^{(\alpha)}=G_n^{(\alpha)}(0)$ are
called the $n$-th Genocchi numbers of order ${\alpha}$. From (1)
and (2), we note that $G_n=G_n^{(1)}$.
 For $d \in \mathbb{N}$ with
$d \equiv 1 \pmod{2}$, let $\chi$ be the  Dirichlet character with
conductor $d$.  It is known that the generalized Genocchi
polynomials attached to $\chi$ are defined by
$$ \left( \dfrac{2t
\sum_{a=0}^{d-1}\chi(a)(-1)^a e^{at}}{e^{dt}+1} \right)
e^{xt}=\sum_{n=0}^\infty G_{n,\chi}(x) \dfrac{t^n}{n!}, \text{
(see [3])}. \eqno(3)$$ In the special case $x=0$, $G_{n, \chi}=
G_{n, \chi}(0)$ are called the $n$-th generalized  Genocchi
numbers attached to $\chi$ (see [3, 4, 5, 6]).

For a real or complex parameter $\alpha$, the generalized
higher-order Genocchi polynomials attached to $\chi$ are also
defined by
$$ \left( \dfrac{2t
\sum_{a=0}^{d-1}\chi(a)(-1)^a e^{at}}{e^{dt}+1} \right)^\alpha
e^{xt}=\sum_{n=0}^\infty G_{n,\chi}^{(\alpha)} (x)
\dfrac{t^n}{n!}, \text{ ( see [13])}. \eqno(4)$$ In the special
case $x=0$, $G_{n, \chi}^{(\alpha)}= G_{n, \chi}^{(\alpha)}(0)$
are called the $n$-th generalized  Genocchi numbers attached to
$\chi$ of order $\alpha$ (see [3, 4, 5, 6, 13, 14, 15]). From (3)
and (4), we derive $G_{n, \chi}= G_{n, \chi}^{(1)}$.

Let us assume that $q \in \mathbb{C} $ with $|q|<1$ as an
indeterminate. Then we use the notation
$$[x]_q=\dfrac{1-q^x}{1-q}.$$
The $q$-factorial is defined by
$$[n]_q !=[n]_q [n-1]_q \cdots [2]_q[1]_q$$ and  the Gaussian
binomial coefficient is also defined by
$${\binom{n}{k}}_q= \dfrac{[n]_q !}{[n-k]_q![k]_q!}=\dfrac{[n]_q [n-1]_q \cdots [n-k+1]_q}{[k]_q!},  \text{ (see [2, 5])}. \eqno(5)$$
Note that $$ \lim_{q \rightarrow
1}\binom{n}{k}_q=\binom{n}{k}=\dfrac{n(n-1)\cdots (n-k+1)}{k!}.$$
It is known that
$$\binom{n+1}{k}_q=\binom{n}{k-1}_q+ q^k \binom{n}{k}_q= q^{n-k} \binom{n}{k-1}_q+\binom{n}{k}_q, \text{ (see [2, 5])}.$$
The $q$-binomial formula are known that
$$(x-y)_q^n=(x-y)(x-qy)\cdots  (x-q^{n-1}y)= \sum_{i=0}^n {\binom{n}{i}}_q q^{{\binom{i}{2}}} (-1)^i x^{n-i} y^i,  \text{ (see [2, 9])},$$
and
$$\dfrac{1}{(x-y)_q^n}=\dfrac{1}{(x-y)(x-qy)\cdots  (x-q^{n-1}y)}
= \sum_{l=0}^\infty  {\binom{n+l-1}{l}}_q   x^{n-l} y^l, \text{
(see [2, 9])}. \eqno(6)$$

There is an unexpected connection with $q$-analysis and quantum
groups, and thus with noncommutative geometry $q$-analysis is a
sort of $q$-deformation of the ordinary analysis. Spherical
functions on quantum groups are $q$-special functions. Recently,
many authors have studied the $q$-extension in various area( see
[1-15]). Govil and Gupta [2] has introduced a new type of
$q$-integrated Meyer-K\"{o}nig-Zeller-Durrmeyer operators and
their results are closely related to study $q$-Bernstein
polynomials and $q$-Genocchi polynomials, which are treated in
this paper. In this paper, we first consider the $q$-extension of
the generating function for the higher-order generalized Genocchi
numbers and polynomials attached to $\chi$. The purpose of this
paper is to present a systemic study of some families of
higher-order generalized $q$-Genocchi numbers and polynomials
attached to $\chi$ by using the generating function of those
numbers and polynomials.

\bigskip
\begin{center} {\bf 2. Generalized $q$-Genocchi numbers and polynomials} \end{center}
\bigskip

For $r \in \mathbb{N}, $ let us consider the $q$-extension of the
generalized Genocchi polynomials of order $r$ attached to $\chi$
as follows:
$$ F_{q, \chi}^{(r)}(t, x)=2^r t^r  \sum_{m_1, \cdots, m_r =0}^{
\infty}  \left(  \prod_{j=1}^r \chi( m_j) \right)
(-1)^{\sum_{j=1}^r m_j} e^{[x+m_1+\cdots+m_r]_q
t}=\sum_{n=0}^\infty G_{n, \chi, q}^{(r)}(x)  \dfrac{t^n}{n!} .
\eqno(7)$$

Note that
$$ \lim_{q \rightarrow 1} F_{q, \chi}^{(r)}(t, x)=\left( \dfrac{2t
\sum_{a=0}^{d-1}\chi(a)(-1)^a e^{at}}{e^{dt}+1} \right)^r
e^{xt}.$$ By (7) and (4), we can see that $ \lim_{q \rightarrow
1}G_{n, \chi, q}^{(r)}(x) =G_{n, \chi}^{(r)}(x).$  From (7), we
note that
$$ \aligned  &   G_{0, \chi,q}^{(r)}(x) =G_{1, \chi, q}^{(r)}(x)=\cdots = G_{r-1, \chi, q}^{(r)}(x)=0,  \text{ and } \\
& \dfrac{ G_{n+r,\chi, q}^{(r)}(x)}{\binom{n+r}{r} r!}=  2^r
 \sum_{m_1, \cdots, m_r =0}^{
\infty}  \left(  \prod_{j=1}^r \chi( m_j) \right)
(-1)^{\sum_{j=1}^r m_j} [x+ m_1+\cdots+m_r]_q^n  .\endaligned
 $$
In the special case  $x=0$, $G_{n, \chi, q}^{(r)}= G_{n, \chi, q
}^{(r)}(0)$ are called the $n$-th generalized  $q$-Genocchi
numbers of order $r$ attached to $\chi$. Therefore, we obtain the
following theorem.

\medskip
{ \bf Theorem 1.} For $r \in \mathbb{N},$ we have
$$  \dfrac{ G_{n+r,\chi, q}^{(r)}}{\binom{n+r}{r} r!}=  2^r
 \sum_{m_1, \cdots, m_r =0}^{
\infty}  \left(  \prod_{i=1}^r \chi( m_i) \right)
(-1)^{\sum_{j=1}^r m_j} [ m_1+\cdots+m_r]_q^n.$$
\medskip

Note that
$$ \aligned  &  2^r
 \sum_{m_1, \cdots, m_r =0}^{
\infty}  \left(  \prod_{i=1}^r \chi( m_i) \right)
(-1)^{\sum_{j=1}^r m_j} [ m_1+\cdots+m_r]_q^n \\
&= \dfrac{2^r}{(1-q)^n}  \sum_{l=0}^n \binom nl (-1)^l \sum_{a_1,
\cdots, a_r =0}^{d-1}  \left(  \prod_{j=1}^r \chi( a_j) \right)
\dfrac{(-q^l)^{\sum_{i=1}^r a_i}}{(1+q^{ld})^r}.
\endaligned
 $$
Thus we obtain the following corollary.

\medskip
{ \bf Corollary 2.} For $r \in \mathbb{N},$ we have
$$ \aligned  &  \dfrac{ G_{n+r,\chi, q}^{(r)}}{\binom{n+r}{r} r!}= \dfrac{2^r}{(1-q)^n}  \sum_{l=0}^n \binom nl (-1)^l \sum_{a_1,
\cdots, a_r =0}^{d-1}  \left(  \prod_{j=1}^r \chi( a_j) \right)
\dfrac{(-q^l)^{\sum_{i=1}^r a_i}}{(1+q^{ld})^r} \\
&= {2^r} \sum_{m=0}^\infty \binom{m+r-1}{m}(-1)^m \sum_{a_1,
\cdots, a_r =0}^{d-1} (-1)^{\sum_{i=1}^r a_i} \left( \prod_{i=1}^r
\chi( a_i) \right) \left[  \sum_{i=1}^r a_i+md \right]_q^n.
\endaligned
 $$

For $h \in \mathbb{Z}$ and $r\in \mathbb{N}$, we also consider the
extended higher-order generalized $(h, q)$-Genocchi polynomials as
follows:

$$ \aligned   F_{q, \chi}^{(h, r)}(t, x) & =2^r t^r  \sum_{m_1, \cdots, m_r =0}^{
\infty}  q^{\sum_{j=1}^r (h-j)m_j} \left(  \prod_{i=1}^r \chi(
m_i) \right) (-1)^{\sum_{j=1}^r m_j} e^{[x+\sum_{j=1}^r m_j]_q t}
\\
& =\sum_{n=0}^\infty G_{n, \chi, q}^{(h, r)}(x)  \dfrac{t^n}{n!} .
\endaligned
\eqno(8)$$

From (8), we note that
$$ \aligned  &   G_{0, \chi,q}^{(h,r)}(x) =G_{1, \chi, q}^{(h,r)}(x)=\cdots = G_{r-1, \chi, q}^{(h,r)}(x)=0,  \text{ and } \\
& \dfrac{ G_{n+r,\chi, q}^{(h,r)}(x)}{\binom{n+r}{r} r!}=  2^r
 \sum_{m_1, \cdots, m_r =0}^{
\infty} q^{\sum_{j=1}^r (h-j)m_j}  \left(  \prod_{i=1}^r \chi(
m_i) \right) (-1)^{\sum_{j=1}^r m_j} [x+ m_1+\cdots+m_r]_q^n \\
&= \dfrac{2^r}{(1-q)^n}  \sum_{l=0}^n \binom nl q^{lx}(-1)^l
\sum_{a_1, \cdots, a_r =0}^{d-1}  \left(  \prod_{j=1}^r \chi( a_j)
\right)q^{\sum_{j=1}^r (h-j)a_j} (-1)^{a_1+\cdots+a_r}
q^{l(a_1+\cdots+a_r)} \\
& \quad \quad \quad \quad  \times \sum_{m_1, \cdots, m_r =0}^{
\infty} (-1)^{m_1+ \cdots + m_r}  q^{d (m_1+\cdots + m_r)+ d( \sum_{j=1}^r (h-j)m_j)} \\
&=\dfrac{2^r}{(1-q)^n}  \sum_{l=0}^n  \dfrac{\binom nl
q^{lx}(-1)^l \sum_{a_1, \cdots, a_r =0}^{d-1}  \left(
\prod_{j=1}^r \chi( a_j) \right)q^{\sum_{j=1}^r (h-j)a_j}
(-q^l)^{\sum_{i=1}^r a_i}}{(-q^{d(h-r+l)}; q)_r}
 ,\endaligned
 $$
where $(-x;q)_r=(1+x)(1+xq) \cdots (1+xq^{r-1}).$

Therefore, we obtain the following theorem.

\medskip
{ \bf Theorem 3.} For $ h \in \mathbb{Z}, r \in \mathbb{N},$ we
have
$$\aligned  & \dfrac{ G_{n+r,\chi, q}^{(h,r)}(x)}{\binom{n+r}{r} r!} =  2^r
 \sum_{m_1, \cdots, m_r =0}^{
\infty} q^{\sum_{j=1}^r (h-j)m_j}  \left(  \prod_{i=1}^r \chi(
m_i) \right) (-1)^{\sum_{j=1}^r m_j} [x+ m_1+\cdots+m_r]_q^n \\
&=\dfrac{2^r}{(1-q)^n}  \sum_{l=0}^n  \dfrac{\binom nl (-q^x)^{l}
\sum_{a_1, \cdots, a_r =0}^{d-1}  \left( \prod_{j=1}^r \chi( a_j)
\right)q^{\sum_{j=1}^r (h-j)a_j} (-q^l)^{\sum_{i=1}^r
a_i}}{(-q^{d(h-r+l)}; q)_r}, \\
& \quad  \text{ and } \quad   G_{0, \chi,q}^{(h,r)}(x) =G_{1,
\chi, q}^{(h,r)}(x)=\cdots = G_{r-1, \chi, q}^{(h,r)}(x)=0.
 \endaligned
 $$
  \medskip

Note that
$$\dfrac{1}{(-q^{d(h-r+l)}; q)_r}=\dfrac{1}{(1+q^{d(h-r+l)})}=\sum_{m=0}^\infty \binom{m+r-1}{m}_q(-1)^m q^{d(h-r+l)m}.\eqno(9)$$
By (9), we see that

$$ \aligned  &  \dfrac{1}{(1-q)^n}  \sum_{l=0}^n  \dfrac{\binom nl (-1)^l   q^{l( x+ \sum_{i=1}^r
a_i)}}{(-q^{d(h-r+l)}; q)_r}  \\
& =\sum_{m=0}^\infty \binom{m+r-1}{m}_q(-1)^m q^{d(h-r)m}
\dfrac{1}{(1-q)^n}  \sum_{l=0}^n  \binom nl (-1)^l   q^{l( x+
\sum_{i=1}^r a_i +dm)} \\
&= \sum_{m=0}^\infty \binom{m+r-1}{m}_q(-1)^m q^{d(h-r)m} [ x+
\sum_{i=1}^r a_i+ dm ]_q^n.
 \endaligned
 \eqno(10)
 $$
By (9) and (10), we obtain the following corollary.

\medskip
{ \bf Corollary 4.} For $ h \in \mathbb{Z}, r \in \mathbb{N},$ we
have
$$\aligned  & \dfrac{ G_{n+r,\chi, q}^{(h,r)}(x)}{\binom{n+r}{r} r!}  \\
&=  2^r
 \sum_{m=0}^\infty \binom{m+r-1}{m}_q(-1)^m q^{d(h-r)m}
\sum_{a_1, \cdots, a_r =0}^{d-1}  \left( \prod_{j=1}^r \chi( a_j)
\right)q^{\sum_{j=1}^r (h-j)a_j} [ x+ \sum_{i=1}^r a_i+ dm ]_q^n.
 \endaligned
 $$
  \medskip

By (8), we can derive the following corollary.

\medskip
{ \bf Corollary 5.} For $ h \in \mathbb{Z}, r , d \in \mathbb{N}$
with $d \equiv 1(\mod 2)$,  we have
 $$\aligned  &  q^{d(h-1)}  \dfrac{ G_{n+r,\chi, q}^{(h,r)}(x+d)}{\binom{n+r}{r} r!} + \dfrac{ G_{n+r,\chi, q}^{(h,r)}(x)}{\binom{n+r}{r} r!}
 = 2\sum_{l=0}^{d-1} \chi(l) (-1)^l \dfrac{ G_{n+r-1,\chi, q}^{(h-1,r-1)}}{\binom{n+r-1}{r-1}
 (r-1)!}, \\
 & \text{  and } \\
 & q^{x}  \dfrac{ G_{n+r,\chi, q}^{(h+1,r)}(x)}{\binom{n+r}{r} r!} = (q-1)  \dfrac{ G_{n+r+1,\chi, q}^{(h,r)}(x)}{\binom{n+r+1}{r} r!}
 + \dfrac{ G_{n+r,\chi, q}^{(h,r)}(x)}{\binom{n+r}{r} r!}.
  \endaligned
 $$
 \medskip

For $h=r$ in Theorem 3, we obtain the following corollary.

\medskip
{ \bf Corollary 6.} For $ r \in \mathbb{N} $,  we have
 $$\aligned  & G_{n+r,\chi, q}^{(r,r)}(x)  \\
&=\dfrac{2^r}{(1-q)^n}  \sum_{l=0}^n \binom nl (-q^x)^{l}
\sum_{a_1, \cdots, a_r =0}^{d-1}  \left( \prod_{j=1}^r \chi( a_j)
\right)  \dfrac{ q^{\sum_{j=1}^r ( (r-j)a_j+ l a_j)} (-1)^{
a_1+\cdots+a_r
}}{(-q^{d l}; q)_r} \\
& = 2^r
 \sum_{m=0}^\infty \binom{m+r-1}{m}_q(-1)^m \sum_{a_1, \cdots, a_r =0}^{d-1}  \left( \prod_{j=1}^r \chi( a_j)
\right)q^{\sum_{j=1}^r (r-j)a_j} [ x+ \sum_{i=1}^r a_i+ dm ]_q^n. \\
& \text{ In particular,} \\
& \dfrac{ G_{n+r,\chi, q^{-1} }^{(r,r)}(r-x)}{\binom{n+r}{r} r!} =
(-1)^n q^{n+ \binom r2} \dfrac{ G_{n+r,\chi,
q}^{(r,r)}(x)}{\binom{n+r}{r} r!}.
 \endaligned
 $$
 \medskip

Let $x=r$ in Corollary 6. Then we have
$$ \dfrac{ G_{n+r,\chi,
q^{-1} }^{(r,r)}}{\binom{n+r}{r} r!} = (-1)^n q^{n+ \binom r2}
\dfrac{ G_{n+r,\chi, q}^{(r,r)}(r)}{\binom{n+r}{r} r!}.
  $$

Let $ w_1, w_2, \cdots, w_r \in \mathbb{Q}_+.$ Then we define
Barnes' type generalized $q$-Genocchi polynomials attached to
$\chi$ as follows:

$$ \aligned   F_{q, \chi}^{(r)}(t, x \mid w_1, w_2, \cdots, w_r ) & =2^r t^r  \sum_{m_1, \cdots, m_r =0}^{
\infty}  \left(  \prod_{i=1}^r \chi( m_i) \right) (-1)^{ m_1+
\cdots + m_r} e^{[x+\sum_{j=1}^r w_j m_j]_q t}
\\
& =\sum_{n=0}^\infty G_{n, \chi, q}^{( r)}(x \mid w_1, w_2,
\cdots, w_r)  \dfrac{t^n}{n!} .
\endaligned
\eqno(11)$$
By (11), we see that
 $$ \dfrac{ G_{n+r,\chi, q}^{(r)} (x \mid w_1,
\cdots, w_r)}{\binom{n+r}{r} r!}=  2^r
 \sum_{m_1, \cdots, m_r =0}^{
\infty}  \left(  \prod_{i=1}^r \chi( m_i) \right)
(-1)^{\sum_{j=1}^r m_j} [ x+\sum_{j=1}^r w_j m_j]_q^n.$$

It is easy to see that

$$ \aligned  & 2^r
 \sum_{m_1, \cdots, m_r =0}^{
\infty}  \left(  \prod_{i=1}^r \chi( m_i) \right) (-1)^{ m_1+
\cdots + m_r} [ x+\sum_{j=1}^r w_j m_j]_q^n \\
& =\dfrac{2^r}{(1-q)^n}  \sum_{l=0}^n  \dfrac{\binom nl (-q^x)^{l}
\sum_{a_1, \cdots, a_r =0}^{d-1}  \left( \prod_{j=1}^r \chi( a_j)
\right) (-1)^{\sum_{j=1}^r a_j} q^{ l \sum_{i=1}^r w_i
a_i}}{(1+q^{dlw_1}) \cdots (1+q^{dlw_r})}  .
\endaligned
$$

Therefore, we obtain the following theorem.

\medskip
{ \bf Theorem 7.} For $ r \in \mathbb{N}, w_1, w_2, \cdots, w_r
\in \mathbb{Q}_+, $ we have
$$\aligned  & \dfrac{ G_{n+r,\chi, q}^{(r)}(x \mid w_1, w_2,
\cdots, w_r)}{\binom{n+r}{r} r!} \\
&  =  2^r
 \sum_{m_1, \cdots, m_r =0}^{
\infty}   \left(  \prod_{i=1}^r \chi(
m_i) \right) (-1)^{\sum_{j=1}^r m_j} [x+ w_1m_1+\cdots+ w_r m_r]_q^n \\
&=\dfrac{2^r}{(1-q)^n}  \sum_{l=0}^n  \dfrac{\binom nl (-q^x)^{l}
\sum_{a_1, \cdots, a_r =0}^{d-1}  \left( \prod_{j=1}^r \chi( a_j)
\right) (-1)^{\sum_{j=1}^r a_j} q^{ l \sum_{i=1}^r w_i
a_i}}{(1+q^{dlw_1}) \cdots (1+q^{dlw_r})}  .
 \endaligned
 $$
  \medskip

\bigskip

\begin{center}{\bf REFERENCES}\end{center}
\begin{enumerate}

\item
{ I. N. Cangul, V. Kurt, H. Ozden, Y. Simsek}, { On the
higher-order $w$-$q$-Genocchi numbers,} { Adv. Stud. Contemp.
Math.}, {19}(2009), 39-57.

\item
{N. K. Govil, V. Gupta,} {Convergence of
$q$-Meyer-K\"{o}nig-Zeller-Durrmeyer operators,} { Adv. Stud.
Contemp. Math.}, {19}(2009), 97-108.

\item
{L.-C. Jang, K.-W. Hwang, Y.-H. Kim,} {A note on $(h, q)$-Genocchi
polynomials and numbers of higher order,} { Adv. Difference.
Equ.}, 2010, Art ID 309480, 6pp.

\item
{L.-C. Jang,} {A  study on the distribution of twisted
$q$-Genocchi polynomials,} { Adv. Stud. Contemp. Math.},
{18}(2009), 181-189.

\item
{T. Kim,}  {  Some identities for the Bernoulli, the Euler and the
Genocchi numbers and polynomials,} { Adv. Stud. Contemp. Math.},
{20}(2010), 23-28.

\item
{T. Kim,} { A note on the $q$-Genocchi numbers and polynomials,}
{ J. Inequal. Appl.}, 2007,  Art ID 71452, 8pp.

\item
{ T. Kim,}  {   On the multiple $q$-Genocchi and Euler numbers},{
Russ. J. Math. Phys.},  { 15}(2008), 481-486.

\item
{T. Kim,} { On the $q$-extension of Euler and Genocchi numbers,} {
J. Math. Anal. Appl.}, {326}(2007) 1458-1465.

\item
{T. Kim,} { Barnes type multiple $q$-zeta function and $q$-Euler
polynomials,} { J. Phys. A: Math. Theor.}, {43}(2010) 255201,
11pp.

\item
{T. Kim,}  {  Note on the Euler $q$-zeta functions,} {J.  Number
Theory}, { 129}(2009), 1798-1804.

\item
{V. Kurt}, {A further symmetric relation on the analogue of the
Apostol-Bernoulli and the analogue of the Apostol-Genocchi
polynomials,} { Appl. Math. Sci.},  { 3}(2009), 2757-2764.

\item
{M. Cenkci, M. Can,  V. Kurt}, {$q$-extension of Genocchi
numbers,} { J. Korean Math. Soc.}, { 43}(2006), 183-198.

\item
{S.-H. Rim, S. J. Lee, E. J. Moon, J. H. Jin}, { On the
$q$-Genocchi numbers and polynomials associated with $q$-zeta
function}, { Proc. Jangjeon Math. Soc.}, 12(2009), 261-267.

\item
{S.-H. Rim, K. H. Park, E. J. Moon}, { On Genocchi numbers and
polynomials}, { Abstr. Appl. Anal.},  Art ID 898471, 7pp.

\item
{C. S. Ryoo}, { Calculating zeros of the twisted Genocchi
polynomials,} { Adv. Stud. Contemp. Math.}, {17}(2008), 147-159.

\end{enumerate}

\end{document}